\documentclass[10pt]{article}
\title{On normal forms for the similarity classes  of matrices and pairs of matrices}
\author{Klaus Bongartz}\date{klausbongartz@online.de}
\usepackage{amsthm} 
\usepackage[utf8]{inputenc}
\newtheorem{definition}{Definition}
\newtheorem{theorem}{Theorem}
\newtheorem{lemma}{Lemma}
\newtheorem{corollary}{Corollary}
\newtheorem{proposition}{Proposition}

\begin{document}
\maketitle
\begin{abstract}We answer two questions posed 1998 in the book `Arnold's problems'. First, over any field $k$ there is  a representative system for the similarity classes of $n\times n$-matrices  which is a finite disjoint union of affine subspaces.  And second, for $n\geq 2$ an analogous statement fails for pairs of $ n\times n$- matrices over any algebraically closed field of characteristic 0. The positive answer rests on a slight modification of the rational normal form whereas the negative answer requires more work.
\end{abstract}

On page 127  of \cite{A} Kontsevich asks whether there is in the space of complex matrices  of a fixed size  a representative system for the similarity classes which is the union of affine subspaces. In the comments he writes on page 613  that this has an easy positive answer  for one matrix which he does not remember but that the analogous question for pairs of matrices is also  interesting and open.

In this note we describe first an easy solution  for one matrix given already 1989 in \cite{B1}, then we give the negative answer  for pairs of  $2 \times 2$-matrices using a result of \cite{B2} from 1996 and finally we extend the answer to all $n\geq 3$ using appropriate vector bundles as in \cite{B3}.

We fix some notations and we recall basic facts about the rational normal form.
Let $k$ be an arbitrary  field. The set of all $n\times n$-matrices with coefficients in $k$ is denoted by $k^{n\times n}$ and the canonical base vectors in $k^{n}$ by $e_{i}$. 

\begin{definition}A rational normal form $R \in  k^{n\times n}$ is a block-diagonal matrix \[
\left [
\begin{array}{cccc}
B(P_{1})&0& ...&0\\
0&B(P_{2})& ... & 0\\
...&...&...&...\\
0&0&...&B(P_{r})
\end{array}
\right ]
\] with companion-matrices $B(P_{1}),B(P_{2}),\ldots ,B(P_{r})$ on the main diagonal, such that all polynomials $P_{i}$ are normed of degree $\geq 1$ with coefficients in $k$ and such that  
 $P_{i+1}$  divides $P_{i}$for each $i$ with  $1\leq i \leq r-1$. We write $R=R(P_{1},P_{2},\ldots ,P_{r})$ for such a matrix.
\end{definition}

\newpage 
There is the following important classical result ( see \cite{B4} for a short proof ).

\begin{theorem}Each matrix $A \in k^{n\times n}$ is similar over $k$ to exactly one rational normal form $R(A)=R(P_{1},P_{2},\ldots ,P_{r})$. Here $R(A)$ as well as an invertible matrix
 $T$ with $T^{-1}AT=R$ are obtained from the coefficients of $A$ by finitely many rational operations the number of which grows polynomially with $n$. 
\end{theorem}              

\begin{definition} Let $p=(p_{1},p_{2}, \ldots ,p_{r})$ be a partition of $n$, i.e. a sequence of natural numbers  $p_{1}\geq p_{2}\geq \ldots \geq p_{r}\geq 1$ adding up to $n$. 
 We denote by $R(p)$ the set of rational normal forms as above with blocks of sizes $p_{1},p_{2}, \ldots ,p_{r}$.
\end{definition}
 
Thus  the disjoint union of the various $R(p)$  is  a representative system for the similarity classes in $k^{n \times n}$. Each $R(p)$ is a closed subset of  $k^{n\times n}$ but not an affine subspace as soon as two different $p_{i}$'s occur.  To remedy this we replace the companion-matrices   by `generalized companion matrices'  as defined below.

\begin{definition} Let $q=(q_{1},q_{2}, \ldots,q_{s})$ be an $s-tuple$ of natural numbers $q_{j} \geq 1$. Given normed polynomials $Q_{j}$ of degree $q_{j}$  the generalized companion-matrix $C=C(Q_{1},Q_{2},\ldots ,Q_{s})$ is the quadratic matrix of size  $m=\sum q_{i}$ with the companion-matrices $B(Q_{j})$ as diagonal blocks ordered in the natural way. Outside of these blocks all entries  are $0$ except that $c_{i+1,i}=1$ holds for all $i<m$.
By $C(q)$ we denote the set of all these generalized companion matrices.
\end{definition}

As an illustration we write down all matrices in $C(q)$ for $q=(2,1,2)$. Here $a,b,c,d,e$ are arbitrary scalars.
 \[
\left [
\begin{array}{cccccccccccc}
0&e&0&0&0\\
1&d&0&0&0\\
0&1&c&0&0\\
0&0&1&0&b\\
0&0&0&1&a\\

\end{array}
\right ]
\].

The following  observation is crucial:
\begin{lemma} The matrix $C=C(Q_{1},Q_{2}, \dots ,Q_{s})$ from above is similar to the companion-matrix of $P=Q_{1}\cdot Q_{2} \cdots Q_{s}$. Furthermore $C(q)$ contains exactly one nilpotent matrix $N$ and the differences $A-N$ with $A \in C(q)$ form  a linear subspace of $k^{m\times m}$ of dimension $m$.
\end{lemma}
Proof:   We make $k^{m}$ to a $k[X]$-module $M$ via  $X\cdot v = Cv$. It is enough to prove  that $e_{1}$  generates $M$ and that $P$ annihilates $e_{1}$.  Namely then the map $1 \mapsto e_{1}$ induces a surjection  $k[X]/(P) \rightarrow M$ between spaces of the same dimension and so an isomorphism. Thus $B(P)$ is similar to $C$.

We proceed by induction on $s$. The case $s=1$ is obvious. To go from $s-1$ to $s$  suppose  $Q_{1}=a_{0}+ a_{1}X + \ldots a_{q_{1}-1}X^{q_{1}-1} + X^{q_{1}}$.  We have $Xe_{1}=e_{2}, Xe_{2}=e_{3},\ldots ,Xe_{q_{1}-1}=e_{q_{1}}$
and $Xe_{q_{1}}= e_{q_{1}+1} -(a_{0}e_{1}+ a_{1}e_{2} + \ldots a_{q_{1}-1}e_{q_{1}} )$. Thus we get  $Q_{1}e_{1}=e_{q_{1}+1}$. 
The action of $X$ on the submodule $U$ of $M$ spanned by the
$e_{j}$ with $j \geq q_{1}+1$ is  given by a generalized companion-matrix with $s-1$ blocks, whence by induction $U$ is generated by $e_{q_{1}+1}$ which is annihilated by   $Q_{2}\cdot Q_{3}\cdots Q_{s}$.  So $e_{1}$ generates $M$ and  it is annihilated by  $P$. 

Now $B(P)$ is nilpotent iff $P=X^{m}$. The rest is obvious.   

\newpage Our aim is to prove:

\begin{theorem} Let $p=(p_{1},p_{2}, \ldots ,p_{r})$ be a partition of $n$. Then there is an affine subspace $A(p)$ of dimension  $p_{1}$ in $k^{n \times n}$ such that the map $A \mapsto R(A)$ induces a bijection between $A(p)$ and $R(p)$.
\end{theorem}
Proof:  First we take a closer look at the partition. Assume that there are $s$ different $p_{i}$'s occurring in $p$. In case  $s=1$ we set $A(p)=R(p)$ which is already an affine subspace. 
 For $s\geq 2$ we number the `jump indices' $i<r$  where $p_{i}\neq  p_{i+1} $ in ascending order  as $j_{1}<j_{2}<\ldots< j_{s-1}$ and we add $j(s)=r,j_{0}=0$. Thus we obtain for  $1\leq i \leq s-1$ non-zero natural numbers $q_{i}=p_{j_{i}}-p_{j_{i}+1}$ to which we add $q_{s}=p_{r}$.
 
  For any  $s$-tuple $(Q_{1},Q_{2}, \ldots ,Q_{s})$ of normed polynomials with $deg Q_{j}=q_{j}$  we consider a block-diagonal matrix 
  $A(Q_{1},Q_{2},\ldots Q_{s})$ of the shape
 \[
\left [
\begin{array}{cccc}
D_{1}&0& ...&0\\
0&D_{2}& ... & 0\\
...&...&...&...\\
0&0&...&D_{r}
\end{array}
\right ]
.\]  
 
Here we set $D_{i}=C(Q_{k},Q_{k+1},\ldots ,Q_{s})$ for  $j_{k-1}< i \leq j_{k}$. This matrix depends only on the tuple of polynomials and on the partition $p$. The set of all these matrices is the wanted affine subspace. Namely each block $D_{i}$ of $A(Q_{1},Q_{2},\ldots,Q_{s})$  is just a `right
lower' sub-block of $D_{1}$ the size of which depends only on $p$ and the possible $D_{1}$'s form an affine subspace. 
As an illustration we write down all matrices in $D(p)$ for the partition $(5,3,2,2)$ of $12$. Here $a,b,c,d,e$ are arbitrary scalars.
 \[
\left [
\begin{array}{cccccccccccc}
0&e&0&0&0&0&0&0&0&0&0&0\\
1&d&0&0&0&0&0&0&0&0&0&0\\
0&1&c&0&0&0&0&0&0&0&0&0\\
0&0&1&0&b&0&0&0&0&0&0&0\\
0&0&0&1&a&0&0&0&0&0&0&0\\
0&0&0&0&0&c&0&0&0&0&0&0\\
0&0&0&0&0&1&0&b&0&0&0&0\\
0&0&0&0&0&0&1&a&0&0&0&0\\
0&0&0&0&0&0&0&0&0&b&0&0\\
0&0&0&0&0&0&0&0&1&a&0&0\\
0&0&0&0&0&0&0&0&0&0&0&b\\
0&0&0&0&0&0&0&0&0&0&1&a\\

\end{array}
\right ]
\]

To prove the theorem we  verify that $A\mapsto R(A)$ is the wanted bijection.
The product $P_{i}=Q_{k}Q_{k+1} \cdots Q_{s}$ is normed of degree  $p_{i}$ and $D_{i}$ is similar to $B(P_{i})$ by lemma 1. Obviously  $P_{i}$ divides $P_{i-1}$ for all $i>1$ and so the rational normal form of the constructed matrix is $R(P_{1},P_{2},\ldots ,P_{r})$ and it belongs to $R(p)$. Reversely, given a rational normal form $R=R(P_{1},P_{2},\ldots ,P_{r})$  one gets the $s-1$   quotients $Q_{i}=P_{j_{i}}/P_{j_{i}+1}$ having degree $q_{i}$ and one adds $Q_{s}=P_{r}$. Then $A(Q_{1},Q_{2}, \ldots Q_{s})$ belongs to $A(p)$  and its rational normal form is $R$. 

In \cite{D} Daugulis has given another construction of the affine subspaces starting from the Jordan normal forms so that his proof is more complicated.
I thank him for bringing Kontsevichs problem to my attention.

 In the following all varieties are  quasi-projective and defined  over an algebraically closed field $k$.
\begin{definition} Let $X$ be a variety endowed  with a morphic action $G \times X \rightarrow X$ of an algebraic group $G$. The action admits normal forms if there  is a constructible subset $N$ that hits each orbit in exactly one point. The elements of $N$ are called normal forms. 
\end{definition}
\begin{lemma}In the situation of the definition $X$ admits a set of normal forms provided there exist finitely many constructible subsets $N_{1},N_{2}, \ldots N_{r}$ such that each orbit hits their union at least once, but each $N_{j}$ at most once. 
\end{lemma}
Proof: By a fundamental theorem of Chevalley the sets  $ G( \bigcup _{j >i}N_{j})$ are  for $i=1,2, \ldots r-1$    constructible as images of constructible sets under the  morphism
$G\times X \rightarrow X$. Thus the sets  $N'_{i}=N_{i}\setminus G( \bigcup _{j >i}N_{j})$ are constructible and their union hits each orbit exactly once.
Namely given $x \in X$ let $j$ be the largest index such that the orbit $Gx$ hits $N_{j}$. Then $Gx \cap N'_{j}$ is just one point,  but $Gx \cap  N'_{i}= \emptyset$ for $i\neq j$.\vspace{0.5cm}

 For the action of $GL_{n}(k)$ on $k^{n \times n}$ by conjugation the finitely many sets  $R(p)$ resp. 
 $A(p)$ are closed subsets and  their unions $R$ resp. $A$ are  both  systems of normal forms.

 We want to show:
\begin{theorem} If the characteristic of $k$ is $0$ the action of $GL_{n}(k )$ on $V(n):=({k}^{n \times n})^{2}$ by simultaneous conjugation
does not admit normal forms for $n\geq 2$. 
\end{theorem}

We need some  observations  essentially taken from section 6.2 of \cite{B2}. 
\begin{proposition} Consider the morphism $F:(sl_{2}(k))^{2} \rightarrow k^{3}$ sending a pair $(A,B)$ of two $2 \times 2$ matrices with trace $0$ to the triple $(det A,tr AB, det B)$.  Let $Y$ be the subset of $k^{3}$ where $g=X_{1}X_{3}(X_{2}^{2} - 4 X_{1}X_{3})$ does not vanish and denote by 
$$ f:X:=F^{-1}Y \rightarrow Y$$ the morphism induced by $F$. Furthermore let $Q$ be the subset of $sl_{2}(k)$ consisting of the pairs  of matrices 
 \[
A= \left [
\begin{array}{cc}
a_{11}&1\\
0&-a_{11}\\

\end{array}
\right ]\] 
 \[
B=\left [
\begin{array}{cc}
b_{11}&0\\
b_{21}&-b_{11}\\

\end{array}
\right ]\] where $a_{11}b_{11}b_{21}(4a_{11}b_{11}+b_{21})\neq 0$. Then we have:
\begin{enumerate}

\item $GL_{2}$ acts on $X$ by conjugation and $F$ is $GL_{2}$ invariant.
\item $Q$ belongs to $ X$.
\item Two matrices $A$ and $B$ with $(A,B) \in X$ have no common eigenvector. 
\item The fibres of $f$ are the orbits. Each orbit hits  $Q$ once  for $char k=2$ and four times  for $char k \neq 2$. 
\item $f$ has no local section, i.e. there is no dense open subset $U$ of $Y$  admitting a morphism  $s:U \rightarrow X$ with $f\circ s = id_{U}$. 

\end{enumerate}
\end{proposition}

Proof: The first assertion is obvious and the second a short calculation.  Now  we consider two matrices $A,B$ in $sl_{2}$ having  a common eigenvector $v$, say $Av=av$ and $Bv=bv$. With respect to an arbitrary basis $(v,w)$ of $k^{2}$ then $A$ resp. $B$ are represented by upper triangular matrices with diagonal entries $a,-a$ resp. $b,-b$ and one gets that $g$ vanishes on $F(A,B)= (-a^{2},2ab,-b^{2})$. 

 Given $y=(x_{1},x_{2},x_{3})$ in $Y$ there are  square roots $a_{11}$ of $-x_{1}$ and $b_{11}$ of $-x_{3}$. Then there is exactly one pair
 $(A,B) \in Q$ mapping to $y$ and one has  $b_{21}=x _{2}-2b_{11}a_{11}$. For $char k \neq 2$ one has two different choices for the square roots leading to four elements of $Q$ in $f^{-1}y$ whereas there is only one for $char k=2$.

 Take now an arbitrary $(A,B) \in f^{-1}y$. For $char k \neq 2$ there are two eigenvectors of $A$ with $Av_{1}=av_{1}$ and  $Av_{2}=-av_{2}$. Similarly we get $Bw_{1}=bw_{1}$ and $Bw_{2}=-bw_{2}$. Since $A$ and $B$ have no common eigenvector $(v_{1},w_{1})$ is a basis of $k^{2}$  and $Aw_{1}=-aw_{1} + xv_{1}$ with $x\neq 0$. With respect to the basis $(xv_{1},w_{1})$  the pair $(A,B)$ is represented by a pair in $Q$. The other three choices $(v_{i},w_{j})$ for bases of $k^{2}$ lead to the three other elements in $Q\cap f^{-1}y$. They are all conjugate since they represent the same endomorphism. For $char k=2$ the reasoning is similar but simpler because $A$ and $B$ have only one eigenvalue.
 
 It remains to prove the last statement which is a bit tricky. So assume $s:U \rightarrow X$ is a local section given by $s(u)=(A(u),B(u))$.
 Fix a point $u_{0}\in U$. Then there is an eigenvecor $v$ of $B(u_{0})$  and this is not an eigenvector of $A(u_{0})$ whence  $(v,A(u_{0})v)$ is a basis of $k^{2}$. The matrix $M(u):=(v,A(u)v)$ is then invertible in an open neighborhood $U'$ of $u_{0}$. Define $s':U' \rightarrow X$  by $s'(u)= (M(u)^{-1}A(u)M(u),M(u)^{-1}B(u)M(u))$. This  is a local section on $U'$ because $f$ is $GL_{2}$-invariant. Replacing $s$ by $s'$ and writing $u=(u_{1},u_{2},u_{3})$ for a point in $U$ we can assume right from the beginning
  that $s(u)=(A(u),B(u))$ with  \[
A(u)= \left [
\begin{array}{cc}
0&-u_{1}\\
1&0\\

\end{array}
\right ]\] 
 \[
B(u)=\left [
\begin{array}{cc}
b_{11}(u)&b_{12}(u)\\
b_{21}(u)&-b_{11}(u)\\

\end{array}
\right ].\]

Restricting even  further we have $b_{ij}(u)=\frac{P_{ij}(u)}{Q(u)}$ for some polynomials $P_{ij},Q$ in $k[X_{1},X_{2}, 
X_{3}]$ with $Q\neq 0$.  Since $s$ is a section and $U$ is  dense in $k^{3}$ we obtain in $k[X_{1},X_{2}, 
X_{3}]$ the two equalities $QX_{2}= -X_{1}P_{21} + P_{12}$ and $Q^{2}X _{3}= -P_{11}^{2}-P_{12}P_{21}$ and so the equality 
$$Q^{2}X_{3}= - P_{11}^{2}-QP_{21}X_{2}- P_{21}^{2}X_{1}.$$  The total degree of the left hand side     is odd and it coincides with the total degree of the right hand side whence $P_{21}\neq 0$. Observing that the total degree of only one term is even one finds  that neither $deg Q > deg P_{21}$ nor $deg Q < deg P_{21}$ is possible. We replace $Q$ and $P_{21}$ by their leading terms $\tilde{Q}$ and $\tilde{P}_{21}$ and we get the equation
$$\tilde{Q}^{2}X_{3}= -\tilde{Q}\tilde{P}_{21}X_{2}- \tilde{P}_{21}^{2}X_{1}$$ of homogeneous polynomials. Now $k[X_{1},X_{2},X_{3}]$ is a unique factorization domain. Any polynomial  $R$ dividing  $\tilde{Q}$ and $\tilde{P}_{21}$ can be cancelled in both. Thus we reduce to the case where they have no common irreducible factor. We infer that $\tilde{Q}=aX_{1}^{t}$ and $\tilde{P}_{21}=bX_{3}^{t}$. This leads to the contradiction
$$a^{2}X_{1}^{2t}X_{3}= -abX_{1}^{t}X_{3}^{t}X_{2}- b^{2}X_{3}^{2t}X_{1}. $$

In the proof of theorem 3  also the next well-known fact is used. 
\begin{proposition} Let $f: X \rightarrow Y$ be an injective dominant morphism between irreducible varieties. For $char k=0$  there is a dense open subset $U$ of $Y$ such that the induced morphism $f^{-1}U \rightarrow U$ is an isomorphism.  
\end{proposition} 

Proof: In characteristic $0$ all field extensions are separable and so we can apply the results presented  in sections 4.6 and 4.7 of \cite{H}. The proofs given there are due to Chevalley. Alternatively one finds a proof in appendix 1 of the nice book \cite{K}.

\begin{corollary} Let $N$ be a system of normal forms for the action of an algebraic group $G$ on an irreducible variety $X$. If there is a morphism $f:X \rightarrow Y$ with the orbits as fibres then $f$ has a local section provided the characteristic is $0$.
\end{corollary}
Proof: As a constructible set $N$ is the union of finitely many locally closed subsets each of which is a finite union of irreducible locally closed subsets. Thus $N$ is the union of finitely many irreducible locally closed subsets $N_{i}$.  Since $Y=f(X)=f(N)$ is  irreducible too there is an index $j$ such that $Y=\overline{f(N_{j})}$ and so the induced morphism $g:N_{j} \rightarrow Y$ is dominant and injective. By proposition 2 there is an open dense subset $U$ of $Y$ such that $h:g^{-1}(U) \rightarrow U$ is an isomorphism. Its inverse induces  the wanted local section.
\vspace{0.5cm}

Now we prove theorem 3 for $n=2$ by contradiction. Assume that $N$ is a set of normal forms in $V(2)$. The subset $X$ of proposition 1  is a $GL_{2}$-stable irreducible subvariety of $V(2)$ and  $X\cap N$ is a set of normal forms for $X$. By proposition 1 we have a morphism
$f:X \rightarrow Y$ with the orbits as the fibres that has no local section. This is a contradiction to the corollary. 
 
For $n\geq 3$ the proof is again by contradiction. We assume that $N$ is a set of normal forms for $V(n)$. Out of $N$ we construct normal forms for $X$ using certain vector bundles as in \cite{B3}. The next picture shows the plan of the proof.\newpage 
The hooked horizontal arrows are  all inclusions of locally closed subsets, $p$ and $q$ are vector bundles, the $s_{i}$'s are a family of morphisms obtained from local sections of $p$ and $q$, the vertical arrow upwards is just an inclusion  and $c$ maps -roughly speaking- a homomorphism to its cokernel. Given a set $N$ of normal forms for $V(n)$  one obtains a set of normal forms for $V(2)'$ out of the $c(s_{i}(N\cap U_{i}))$.

\vspace{0.8cm}

\setlength{\unitlength}{0.8cm}

\begin{picture}(11,4)
\put(5.4,4){$\tilde{Z}\hookleftarrow$}\put(5,3.1){q}
\put(5,2){$\tilde{W}(n)' $}\put(7,0){$U_{i}$}\put(7.2,0.6){\vector(0,1){3.3}}
\put(3,2){$\tilde{W}(n)\hookleftarrow$}\put(3,1.1){p}
\put(1,0){$V(n) \hookleftarrow $}
\put(3,0){$W(n)\hookleftarrow$}
\put(5,0){$W(n)'\hookleftarrow$}
\put(7,4){$Z$}\put(8,4.4){c}\put(7.6,2.2){$s_{i}$}
\put(9,4){$V(2)' \hookrightarrow$}\put(9.2,2){$X$}\put(9.4,2.5){\vector(0,1){1.3}}
\put(11,4){$V(2)$}
\put(5.6,3.8){\vector(0,-1){1.2}}
\put(7.5,4.15){\vector(1,0){1.2}}
\put(3.6,1.8){\vector(0,-1){1.2}}
\put(5.6,1.8){\vector(0,-1){1.2}}
\end{picture}

 The two vector bundles we consider are special cases of the following construction:
\begin{lemma} Let $f:X \rightarrow k^{p\times q}$ be a morphism of varieties. Then the set $X(r)$ of points $x$ with $f(x)$ of rank $r$ is a locally closed subvariety of $X$ for any natural number. Moreover, the closed subset $V=\{(x,v)\mid f(x)v=0\}$ is a subbundle of rank $q-r$ of the trivial bundle $X(r)\times k^{q}$.
\end{lemma}
Proof: The rank of $f(x)$ is $r$ iff the determinants of all quadratic submatrices with $r+1$ columns vanish and  the determinant of at least one $r\times r$-submatrix not. This shows that $X(r)$ is locally closed. For notational reasons we consider  the case where the left upper $r\times r$ submatrix $A_{11}(x)$ of $f(x)$ is invertible. We write $f(x)$ as a block-matrix 

 \[\left [
\begin{array}{cc}
A_{11}(x)&A_{12}(x)\\
A_{21}(x)&A_{22}(x)\\

\end{array}
\right ]\] and $v$ as a corresponding block-column with $v_{1}\in k^{r},v_{2}\in k^{q-r}$.  Then $v$ satisfies $f(x)v=0$ iff $v_{1}=-A_{11}(x)^{-1} A_{12}(x)v_{2}$. Thus on the open set $U$ where $det\,A_{11}(x) \neq 0$  the map $(u,v_{2}) \mapsto (u,v)$ is an isomorphism between $U\times k^{q-r}$ and $p^{-1}(U)$.  \vspace{0.5cm}

Any point $m=(m_{1},m_{2}) \in V(n)$ corresponds to the  $k\langle X_{1},X_{2}\rangle$ - module $M$ over the free associative algebra $A$ in two variables with underlying vectorspace $k^{n}$ where $X_{i}$ acts via $m_{i}$. Two points $m$ and $m'$ are in the same orbit iff the corresponding modules $M$ and $M'$  are isomorphic. For $m\in V(n)$ and $m'\in V(n')$ a homomorphism from $M$ to $M'$ is just a matrix
$f \in k^{n'\times n}$ with $fm=m'f$ i.e.  $fm_{i}=m'_{i}f$ holds for $i=1,2$. Thus $f$ is a solution of a linear equation $B(m,m')f=0$ where the coefficients of $B(m,m') \in k^{2nn' \times nn'}$ depend polynomially on the entries of $m$ and $m'$. It is obvious that $[M,M']:= dim\, Hom_{A}(M,M')= dim \,ker B(m,m')=
nn'- rank \,B(m,m')$.  Thus  for any $M$ and for any natural number $t$ the set of all matrices $m'$ with $[M,M']\geq t$ resp. $[M,M']\leq t$ resp. $[M,M']=t$ is closed resp. open resp. locally closed. An analogous statement holds for  $m$ and  $M'$.

Next we fix one point $s=(s_{1},s_{2})$ in $V(n-2)$ where $s(1)$ is a diagonal matrix having $1,2,\dots ,n-2$ on the diagonal and $s_{2}$ is the permutation matrix to the $n-2$ cycle. Then $S$ is simple. By the above remarks the set
$W(n)=\{m\mid \,[S,M]=1=[M,S] \}$
is a locally closed  subvariety of $V(n)$ and  $\tilde{W}(n)=\{(m,f,g)\mid  m \in W(n),fs=mf,  gm=sg\}$ is  a vector bundle  of rank $2$ over $W(n)$. 
The projection $p$ is open and $GL_{n}$ equivariant for the action $h\ast (m,f,g)= (h^{-1}mh,h^{-1}f,gh)$. Thus $p$ maps the open $GL_{n}$-stable subset
$\tilde{W}(n)'=\{(m,f,g)\mid g\circ f\neq 0\}$  to an open $GL_{n}$-stable subset denoted by $W(n)'$.

For a triple $(m,f,g) \in \tilde{W}(n)'$ one gets by Schurs lemma 
 $g\circ f=a \cdot id_{S}$ for some non-zero scalar $a$, whence $e=a^{-1}f \circ g$ is an idempotent endomorphism of $M$ and so $M$ is the direct sum of $im\,e = im \,f \simeq S$ and $ker\,e = ker\,g$. Over $\tilde{W}(n)'$ there is the vector bundle 
 $$\tilde{Z}=\{(m,f,g,y_{1},y_{2})\mid (m,f,g) \in \tilde{W}(n)', y_{i} \in k^{n}, gy_{1}=gy_{2}=0 \}$$
 which contains the open subset $Z$ where $(y_{1},y_{2})$ is a basis of  $ker\,g$.   Any point $(m,f,g,y_{1},y_{2})\in Z$ produces a basis 
 $h=(fe_{1},fe_{2},\dots ,fe_
 {n-2},y_{1},y_{2})$ of $k^{n}$. The action of $m$ with respect to this basis is given by 
 
  \[
h^{-1}mh=s\oplus t=\left [
\begin{array}{cc}
s&0\\
0&t\\

\end{array}
\right ].\] Here $M\simeq S\oplus T$ and $[S,M]=[M,S]=1$ imply  $[S,T]=[T,S]=0$ and so $s\oplus t$ is similar to $s\oplus t'$ iff $t$ and $t'$ are similar.
 
  We define $c:Z \rightarrow V(2)$ by 
 $c(m,f,g,y_{1},y_{2})=(t_{1},t_{2})=t.$ The image $V(2)'$ of  $c$ is the open subset of $V(2)$ of all $t$ with $[S,T]=[T,S]=0$. The set $X$ defined in proposition 1 consists of simple modules given by matrices with trace $0$ whence not isomorphic to $S$. Thus $X$ is for all $n\geq 3$ contained in $V(2)'$.

 Over any point $m$ in $W(n)'$ we choose  first a point $(m,f,g) \in \tilde{W}(n)'$ and then a point $(m,f,g,y_{1},y_{2}) \in Z$. Then there are local sections  $\sigma:U \rightarrow  \tilde{W}(n)'$ and $\sigma':U' \rightarrow Z$ with open sets $U,U'$, $m \in U$, $(m,f,g) \in U'$ and $\sigma(m)=(m,f,g), \sigma'(m,f,g)=(m,f,g,y_{1},y_{2})$. Then $U\cap \sigma^{-1}V$ is an open neigborhood $U_{m}$ of $m$ and we have a morphism $s_{m}:U_{m} \rightarrow Z$. Since $W(n)'$ is quasi-compact it is covered by finitely many of these $U_{m}$. Changing notations we have for $i=1,2,\ldots r$  morphisms $s_{i}:U_{i} \rightarrow Z$ of the above type so that the open sets $U_{i}$ cover $W(n)'$.

 After all these präparations we  construct now  normal forms for $V(2)'$ out of a set $N$ of normal forms for $V(n)$.
  Since $W(n)'$  is constructible and $GL_{n}$-stable  $N\cap W(n)'$ is a set of normal forms for $W(n)'$. We prove that the sets $N_{i}:=c(s_{i}(N\cap U_{i}))$ satisfy the conditions of lemma 2.
   By a theorem of  Chevalley  all $N_{i}$'s are constructible.  Each $m \in N\cap U_{i}$ is similar to   \[
\left [
\begin{array}{cc}
s&0\\
0&c(s_{i}(m))\\

\end{array}
\right ].\]  Thus  $c(s_{i}(m))$ and  $c(s_{i}(m'))$ are not conjugate under $GL_{2}$ for $m,m' \in N\cap U_{i}$  with $m\neq m'$. Finally for any $t \in V(2)'$ the direct sum $s\oplus t$ is in $W(n)'$ and so it is similar under $GL_{n}$ to an element $m$ in some $N\cap U_{i}$ which is similar to $s\oplus c(s_{i}(m))$. Thus $t$ is similar to an element in $N_{i}$.\vspace{0.6cm}

In a sense this article consists in a recycling of some material contained in old papers of mine. So I make some comments on these  articles and their actual intentions. The study of finite dimensional modules over an associative algebra has categorical and geometrical aspects  with many interactions between them. The main result of \cite{B1} is that the set $S(p)$ of matrices with rational normal form in $R(p)$ is a smooth rational subvariety and the  map $A \mapsto R(A)$ is a smooth morphism which is the geometric quotient. To prove  this purely geometric statement one has to understand better some categorical properties of rational normal forms. Over the complex numbers a theorem similar to mine but without the rationality  was already proven 1978 by Peterson in his Harvard thesis using analytical methods. His geometric quotient  is  an affine subspace  different from $A(p)$ which I introduced  then at the end  of \cite{B1}. 

The study of similarity classes of matrices is just a special  case of the following construction.
Given a finitely generated associative algebra $A =k\langle X_{1},X_{2},\dots X_{n} \rangle /I$ as a quotient of the free associative algebra one has for any $d$ the affine $GL_{d}$-variety $Mod^{d} A$ of $d$-dimensional $A$-modules which consists in $n$-tuples $(M_{1},M_{2},\ldots M_{n})$ of $d\times d$ matrices that satisfy all the relations imposed by $I$. The orbits correspond to the isomorphism classes of $d$-dimensional $A$-modules.  For $d=1$ the varieties $Mod^{d}(A)$ are just the usual affine varieties which are already a set of normal forms. In general the orbits and their closures are difficult to study.
 For the path algebra $A$  of a Dynkin or Euclidean quiver the orbits in $Mod^{d}(A)$ are known  since 1973 by work of Gabriel, Nazarova and Donovan/Freislich \cite{G,N,DF} and  I determined their closures in 1995 \cite{B3}. The first four sections in   \cite{B2} describe  these results shortly. Section 5 gives an introduction into geometric quotients mostly without proofs. Section 6 starts with  a thorough study  of $Mod^{2}k\langle X,Y \rangle$ i.e. of pairs of $2\times 2$- matrices which is more precise than the one in \cite{F}. It  includes  without proof  statement 5 of  proposition 1 which is essential for the proof of theorem 3. The most difficult and interesting result in section 6 is that the set of simple modules in $Mod^{n}k\langle X,Y \rangle$ is the largest open set admitting a geometric quotient. In this proof one stumbles over  an open subset  of $Mod^{3}k\langle X,Y \rangle$ with a free $PGL_{3}$-action but without a geometric quotient. For proper actions two examples of this phenomenon are given by Mumford  on page 83ff.  of \cite{M}. These examples  rely on work of  Nagata and Hironaka. 

Finally, my last article \cite{B4} is only included for publicity. The proof of theorem 1 given there proceeds by induction in three steps the most difficult of which is the first one. The proof  is elementary, constructive and easy to understand for beginners.  Another nice proof using my first step and a duality trick is given by Geck in \cite{G}. He has also written a program.

\end{document}